\documentclass[11pt]{article}

\usepackage[T2A]{fontenc}
\usepackage[cp1251]{inputenc}
\usepackage{amsthm}
\usepackage{amsfonts}
\usepackage{eufrak}
\usepackage{amssymb}
\usepackage{amsmath}
\usepackage{cite}

\makeatletter

\renewcommand{\@seccntformat}[1]{\csname the#1\endcsname. }

\makeatother
\begin{document}
\newtheoremstyle{mytheorem}
  {\topsep}   
  {\topsep}   
  {\itshape}  
  {}       
  {\bfseries} 
  {.}         
  {5pt plus 1pt minus 1pt} 
  {}          
\newtheoremstyle{myremark}
  {\topsep}   
  {\topsep}   
  {\upshape}  
  {}       
  {\bfseries} 
  {.}         
  {5pt plus 0pt minus 1pt} 
  {}          
\theoremstyle{mytheorem}
\newtheorem{theorem}{Theorem}[section]
 \newtheorem{theorema}{Theorem}
 \newtheorem*{heyde1*}{Theorem A}
 \newtheorem*{heyde2*}{Theorem B}
 \newtheorem*{heyde3*}{Theorem C}
 \newtheorem*{heyde4*}{Theorem D}
\newtheorem{proposition}[theorem]{Proposition}
 \newtheorem{lemma}[theorem]{Lemma}
\newtheorem{corollary}[theorem]{Corollary}
\newtheorem{definition}{Definition}[section]
\theoremstyle{myremark}
\newtheorem{remark}[theorem]{Remark}
 \centerline{\textbf{Generalized P\'olya's  theorem on
connected  locally}}

\centerline{\textbf{compact Abelian groups of dimension 1}}

\bigskip

\centerline{\textbf{G.M. Feldman}}

 \bigskip

\noindent {\bf  Abstract.}
According to the generalized
P\'olya theorem, the Gaussian distribution on the real line is
characterized by the property of equidistribution of a   monomial
and a linear  form of independent identically distributed random variables.
We give a complete description  of
 \text{\boldmath $a$}-adic solenoids
 for which an analog  of this theorem is true.
 The proof of the main theorem is reduced to solving some functional equation in the class of continuous positive definite functions on the character group of an \text{\boldmath $a$}-adic solenoid.


\bigskip


\bigskip

\section{Introduction}

According to the classical characterization theorems of mathematical statistics, the Gaussian distribution on the real line is characterized by the independence of two linear forms of independent random variables (the Skitovich--Darmois theorem), as well as by the symmetry of the conditional distribution of one linear form given another   (the Heyde  theorem). A number of papers have been devoted to generalizing of these theorems to various algebraic structures,
 in particular, to locally compact Abelian groups (see e.g.  \cite{Fe-SD-2003, Fe4, Fe20bb, FG, {FG1}, MiFe5, {Fe6}, JOTP2015, Solenoid, M}, and also   \cite[Chapters IV--VI]{Fe5}).
In so doing, the coefficients of   linear forms are topological automorphisms of the group. In this note, we consider a group analog  of the generalized P\'olya theorem, in which the Gaussian distribution on the real line is characterized by the property of equidistribution of a monomial and a linear form of independent identically distributed random variables. We study the generalized P\'olya theorem   on  \text{\boldmath $a$}-adic
solenoids. This is an   important class of locally compact Abelian groups, since each connected locally compact Abelian group of dimension 1 is topologically isomorphic to either the additive group of real numbers, or the circle group, or an \text{\boldmath $a$}-adic
solenoid.
 Furthermore, the coefficients of a monomial and a linear form, as in the group analogs of the Skitovich--Darmois and Heyde theorems, are topological automorphisms of the \text{\boldmath $a$}-adic
solenoid. As far as we know, such problem for groups has not been considered earlier. Note also that in the list of unsolved problems in \cite[Ch. 14]{KaLiRa} the problem of constructing a theory of equidistributed   linear forms on algebraic structures was formulated. A   step in this direction is this note.

Recall the generalized P\'olya theorem.

{\bf Theorem A} (\!\!\cite[\S 13.7]{KaLiRa}).  {\it
Let $\xi_j$, $j=1, 2, \dots, n$, $n\ge 2$, be independent identically
distributed random variables. Let $\alpha_j$ be nonzero real numbers such  that
\begin{equation}\label{41}
 \alpha_1^2+\dots+\alpha_n^2=1.
 \end{equation}
If    $\xi_j$ and     $\alpha_1\xi_1+\dots+\alpha_n\xi_n$ are
identically distributed, then $\xi_j$ are Gaussian.}

If we denote by $f(s)$ the characteristic function of $\xi_j$ it is easy to see that the generalized P\'olya theorem is equivalent to the statement that all solutions of the equation
$$f(s)=f(\alpha_1s)\cdots f(\alpha_ns),
\ \ s\in \mathbb{R},$$
in the class of normed continuous positive definite functions are of the form $f(s)=\exp\{-\sigma s^2+i\beta s\}$, where $\sigma\ge 0$, $\beta\in\mathbb{R}$.

The first theorem of this kind, namely  when $n=2$  and
$\alpha_1=\alpha_2={1\over\sqrt{2}}$, was proved by G. P\'olya in 1923.
Subsequently, the characterization of the Gaussian distribution on the real
line by   the property of  equidistribution  of two linear    forms of
independent identically distributed random variables were studied by
I. Marcinkiewicz, Yu. V. Linnik, A.A. Zinger and others
(see \cite[Chapter 2]{KaLiRa}).

Let $X$ be a second countable locally compact Abelian group. We will
consider only such groups, without mentioning it specifically.
Let ${\rm Aut}(X)$ be the group
of topological automorphisms of the group $X$, and  $I$  be the
identity automorphism of a group. Denote by $Y$ the character
group of the group $X$, and by  $(x,y)$ the value of a character
$y \in Y$ at an element $x \in X$. For a closed subgroup $K$ of $X$,
denote by   $A(Y, K) = \{y \in Y: (x, y) = 1$ for all $x \in K \}$
its annihilator. If $\alpha\in{\rm Aut}(X)$, then
define the adjoint automorphism $\tilde\alpha\in{\rm Aut}(Y)$
by the formula $(x, \tilde\alpha y) = (\alpha x, y)$ for
all $x \in X$, $y \in Y$. Denote by   $\mathbb{R}$ the group of real
numbers, by
 $\mathbb{Z}$ the group of integers and by $\mathbb{T}$ the circle group. Let  $n$ be an integer, $n\ne 0$. Denote by $f_n:X \mapsto X$ an  endomorphism
of the group $X$
 defined by the formula  $f_nx=nx$, $x\in X$.

Let ${\rm M}^1(X)$  be the convolution semigroup of probability
distributions on the group $X$. For  $\mu\in{\rm M}^1(X)$ denote
by $$\hat\mu(y) =
\int_{X}(x, y)d \mu(x), \ \ y\in Y,$$  the characteristic function (Fourier transform)
of the distribution $\mu$, and define the distribution $\bar
\mu \in {\rm M}^1(X)$ by the rule $\bar \mu(B) = \mu(-B)$ for all Borel
subsets  $B$ in $X$.
We have $\hat{\bar{\mu}}(y)=\overline{\hat\mu(y)}$.

A distribution  $\gamma\in {\rm M}^1(X)$  is called Gaussian
(see \cite[Chapter IV, \S 6]{Pa})
if its characteristic function can be represented in the form
$$
\hat\gamma(y)= (x,y)\exp\{-\varphi(y)\}, \ \  y\in Y,
$$
where $x \in X$, and $\varphi(y)$ is a continuous non-negative function
on the group $Y$
 satisfying the equation
 $$
\varphi(u + v) + \varphi(u
- v) = 2[\varphi(u) + \varphi(v)], \ \ u,  v \in
Y.
$$
Note that, in particular, the degenerate
distributions are Gaussian.
Denote by $\Gamma(X)$ the set of Gaussian distributions on
    $X$.

Let $K$ be a compact subgroup of the group $X$. Denote by $m_K$
the Haar distribution on $K$. The
characteristic function of the   distribution $m_K$ is of the
 form
\begin{equation}\label{11a}
\hat m_K(y)=
\begin{cases}
1, & \text{\ if\ }\   y\in A(Y, K),
\\  0, & \text{\ if\ }\ y\not\in
A(Y, K).
\end{cases}
\end{equation}

\section{Main theorem}

Let \text{\boldmath $a$}=$(a_0, a_1,a_2,\dots)$, where all $a_j \in {\mathbb{Z}}$, $a_j > 1$, and let $\Delta_{{\text{\boldmath $a$}}}$ be the
group of  ${\text{\boldmath $a$}}$-adic integers (\!\!\cite[(10.2)]{HeRo1}).  As a
set $\Delta_{{\text{\boldmath $a$}}}$ coincides
with the Cartesian product $\mathop{\mbox{\rm\bf
P}}\limits_{n=0}^\infty\{0,1,\dots ,a_n-1\}$.
Put $\mathbf{u}=(1, 0,\dots,0,\dots)\in \Delta_{{\text{\boldmath $a$}}}$.
Consider the group
$\mathbb{R}\times\Delta_{{\text{\boldmath $a$}}}$ and denote by
 $B$ a subgroup of the group
$\mathbb{R}\times\Delta_{{\text{\boldmath $a$}}}$ of the form
$B=\{(n,n\mathbf{u})\}_{n=-\infty}^{\infty}$. The
factor-group $\Sigma_{{\text{\boldmath $a$}}}
=(\mathbb{R}\times\Delta_{{\text{\boldmath $a$}}})/B$ is called
  an ${\text{\boldmath $a$}}$-{adic solenoid}.   The
group $\Sigma_{{\text{\boldmath $a$}}}$ is  compact, connected and has
dimension 1  (\!\cite[(10.12), (10.13),
(24.28)]{HeRo1}). The character group of the group
$\Sigma_{{\text{\boldmath $a$}}}$ is topologically isomorphic to a discrete additive
group $
 H_{\text{\boldmath $a$}}$ of the rational numbers of the form
 \begin{equation}\label{04.05.1}
H_{\text{\boldmath $a$}}=
\left\{{m \over a_0a_1 \cdots a_n} : \ n = 0, 1,\dots; \ m
\in {\mathbb{Z}} \right\}
\end{equation}
(\!\!\cite[(25.3)]{HeRo1}). It is convenient for us to assume that
$H_{\text{\boldmath $a$}}$ is the character group of the group $\Sigma_\text{\boldmath $a$}$. Each topological automorphism $\alpha\in {\rm
Aut}(\Sigma_\text{\boldmath $a$})$
 is of the following form
$\alpha = f_p f_q^{-1}$ for some mutually prime $p$ and $q$,
where $f_p, f_q \in {\rm
Aut}(\Sigma_\text{\boldmath $a$})$. We will identify $\alpha
=  f_p f_q^{-1}$ with the rational number ${p\over q}$. Note
that $\tilde f_n=f_n$. Hence, if $\alpha = f_p f_q^{-1}\in {\rm
Aut}(\Sigma_\text{\boldmath $a$})$, then $\tilde\alpha = f_p f_q^{-1}\in {\rm
Aut}(H_\text{\boldmath $a$})$.

The main result of this note is the   following theorem.

\begin{theorem}\label{th1} Let $X=\Sigma_{{\text{\boldmath $a$}}}$ be
an {{\text{\boldmath $a$}}}-adic solenoid satisfying the condition:

$(i)$ There is a unique prime number $p$ such that the group $X$
contains no elements of order $p$.

Let $\alpha_j$, $j=1, 2, \dots, n$, $n\ge 2$, be topological automorphisms
of the group   $X$ such that
 \begin{equation}\label{1}
  \alpha_1^2+\dots+\alpha_n^2=I.
 \end{equation}
Let $\xi_j$ be independent identically distributed random variables with
values in    $X$ and distribution  $\mu$. If   $\xi_j$ and     $\alpha_1\xi_1+\dots+\alpha_n\xi_n$ are identically distributed,
then   $\mu=\gamma*m_K$, where  $\gamma\in \Gamma(X)$, and $K$ is a
compact subgroup of the group  $X$. Moreover, $f_p(K)=K$.
\end{theorem}
We note that if $\mu=\gamma*m_K$, where  $\gamma\in \Gamma(X)$
and $K$ is a compact subgroup of the group  $X$, then $\mu$ is invariant
with respect to $K$, and $\mu$ induces on the factor-group   $X/K$ under the natural mapping $X\rightarrow X/K$ a
Gaussian distribution. The condition (\ref{1}) for {{\text{\boldmath $a$}}}-adic
solenoids is an analog  of the condition (\ref{41}) for the real line. Thus,
Theorem \ref{th1}   can be considered   as an analog
for  {{\text{\boldmath $a$}}}-adic solenoids satisfying the condition
$(i)$ of the generalized P\'olya  theorem. Observe that group analogs of the generalized P\'olya  theorem for locally compact
Abelian groups in the case  of equidistribution of a   monomial
and a linear  form with integer coefficients   were studied in \cite[\S 11]{Fe9}.

\begin{corollary}\label{co1} Let $X=\Sigma_{{\text{\boldmath $a$}}}$
be an {{\text{\boldmath $a$}}}-adic solenoid satisfying the  condition
$(i)$ of Theorem $\ref{th1}$. Let $\alpha_j$, $j=1, 2, \dots, n$,
$n\ge 2$, be   topological automorphisms of the group   $X$  satisfying
the condition $(\ref{1})$.
Let $\xi_j$ be independent identically distributed random variables
with values in    $X$ and distribution  $\mu$ with a non-vanishing
characteristic function.
If   $\xi_j$ and     $\alpha_1\xi_1+\dots+\alpha_n\xi_n$
are identically distributed, then    $\mu\in \Gamma(X)$.
\end{corollary}
 \begin{remark}\label{r2} Let \text{\boldmath $a$}=$(a_0, a_1,a_2,\dots)$, and let $\Sigma_{{\text{\boldmath $a$}}}$ be the corresponding {{\text{\boldmath $a$}}}-adic
solenoid.

It is easy to verify
that the following statements are equivalent:

$(i)$ There is a unique prime number $p$ such that   $\Sigma_{{\text{\boldmath $a$}}}$ contains no elements of order $p$;

$(ii)$   There is a unique prime number $p$  such that infinite number
of the numbers   $a_j$ are divided by $p$.

An example of  an  {{\text{\boldmath $a$}}}-adic solenoid
$\Sigma_{{\text{\boldmath $a$}}}$ that satisfies the condition
$(i)$ of Theorem \ref{th1} is an  {{\text{\boldmath $a$}}}-adic
solenoid $\Sigma_{{\text{\boldmath $a$}}}$, where \text{\boldmath $a$}=$(p, p, p,\dots)$. Its character
group is topologically isomorphic to a discrete additive
group $F$ of the rational numbers of the form
\begin{equation}\label{26_09_1}
F=\left\{{m\over p^n}:   \ n = 0, 1,\dots; \ m
\in {\mathbb{Z}}\right\}.
 \end{equation}

 It is easy to see that
the following statements are equivalent:

$(\alpha)$ There are  topological automorphisms   $\alpha_j$, $j=1, 2, \dots, n$,
$n\ge 2$,  of  $\Sigma_{{\text{\boldmath $a$}}}$ satisfying the condition
$(\ref{1})$;

$(\beta)$ There is a  prime number
$p$ such that    $\Sigma_{{\text{\boldmath $a$}}}$ contains no elements
of order  $p$.

$(\gamma)$ There is a  prime number
$p$ such that infinite number of the numbers   $a_j$ are divided by $p$.

To prove it,   note that  there is an alternative:  either $f_p\in {\rm Aut}(\Sigma_{{\text{\boldmath $a$}}})$ for a  prime number
$p$ or ${\rm Aut}(\Sigma_{{\text{\boldmath $a$}}})=\{\pm I\}$.
\end{remark}
To prove Theorem \ref{th1} we need two lemmas.
\begin{lemma}\label{le2} Let  $X$ be a locally compact Abelian group,
and let $\alpha_j$, $j=1, 2, \dots, n$, $n\ge 2$, be
topological automorphisms of the group  $X$. Let $\xi_j$ be independent
identically distributed random variables with values in    $X$ and
distribution  $\mu$. Then    $\xi_j$ and   $\alpha_1\xi_1+\dots+\alpha_n\xi_n$ are identically
distributed if and only if the characteristic function
$\hat\mu(y)$ satisfies the equation
\begin{equation}\label{2}
\hat\mu(y)=\hat\mu(\tilde\alpha_1y)\cdots\hat\mu(\tilde\alpha_ny),
\ \ y\in Y.
 \end{equation}
\end{lemma}
The proof of the lemma is the same as in the classical case and follows directly from the fact that
$\hat\mu(y)={\bf E}[(\xi_j, y)]$ for all $y\in Y$, and the independence
of the random variables  $\xi_j$.

The following lemma can be considered as an  analog for locally
compact Abelian groups of the well-known Cram\'er theorem on
decomposition of the Gaussian distribution on the real line
(\!\!\cite{Fe1},  see also \cite[Theorem 4.6]{Fe5}).
\begin{lemma}\label{le1}     Let  $X$ be a locally compact Abelian
group containing no subgroup topologically isomorphic to the circle
group $\mathbb{T}.$  Let
 $\gamma \in \Gamma(X)$. If
  $\gamma=\lambda_1*\lambda_2$,
 where $\lambda_j\in {\rm M}^1(X)$,
 then $\lambda_j \in \Gamma(X)$, $j=1, 2$.
\end{lemma}
{\bf Proof of Theorem \ref{th1}.} In order not to complicate the notation we will
identify $Y$ with $H_{\text{\boldmath $a$}}$ and  consider  $Y$ as
a subset of the real line $\mathbb{R}$.
 It follows from the definition of the Gaussian distribution that
$\gamma\in \Gamma(X)$ if and only if its characteristic function can
be represented in the form
$$
\hat\gamma(y)= (x,y)\exp\{-\sigma y^2\}, \ \  y\in Y,
$$
where $x \in X$, $\sigma\ge 0$.

By Lemma \ref{le2},  the characteristic function  $\hat\mu(y)$
satisfies equation (\ref{2}).
First, we will prove the theorem assuming that
$\hat\mu(y)\ge 0$ for all $y\in Y$. This implies that  $\hat\mu(-y)=\hat\mu(y)$ for all $y\in Y$.

Since $X$ is a connected group,   $f_n$ is an epimorphism for each nonzero integer $n$
(\!\!\cite[(24.25)]{HeRo1}).
Since $X$ contains no elements of order
$p$, we have that $f_p$ is a monomorphism. In view of the fact that $f_p$ is an epimorphism, $f_p\in {\rm Aut}(X)$,  and
for this reason  $f_p\in {\rm Aut}(Y)$. It means that both $X$ and $Y$ are
groups with unique division by $p$. Note that in view of the condition $(i)$ of
the theorem, if $q$ is a prime number and $q\ne p$, then $f_q\notin {\rm Aut}(X)$. This implies that each topological automorphism  of the
group $X$ is of the form $\pm f^{\pm 1}_{p^m}$, where $m$ is
a non-negative integer. Taking into account (\ref{1}), we conclude that
each topological automorphism $\alpha_j$ in the
theorem is $\alpha_j=\pm f^{-1}_{p^{m_j}}$, where $m_j$ is
a natural number, i.e. $\alpha_j$ is a multiplication by a number  $\pm p^{-m_j}$.  Put $f(y)=\hat\mu(y)$.
If necessary, changing the numbering of the random variables $\xi_j$ and
taking into account that $f(-y)=f(y)$, we can rewrite equation  (\ref{2})
in the form
 \begin{equation}\label{3}
f(y)=(f(p^{-1}y))^{k_1}\cdots (f(p^{-l}y))^{k_l}, \ \ y\in Y,
 \end{equation}
 where $k_j\ge 0$, $j=1, 2, \dots, l$, $k_1+\dots+k_l=n$.   It follows
from  (\ref{1}) that
  \begin{equation}\label{4}
\sum_{j=1}^l{k_j\over p^{2j}}=1.
 \end{equation}
Assume that for each $j$ the inequality  $k_j\le p^j$ is valid. Then
we have
 $$
 \sum_{j=1}^l{k_j\over p^{2j}}\le \sum_{j=1}^l{1\over p^{j}}<1,
 $$
 contrary to (\ref{4}). Hence, there is $j_0$ such that
 \begin{equation}\label{21}
k_{j_0}>p^{j_0}.
 \end{equation}
Taking this into account and the fact that $0\le f(y)\le 1$, (7) implies the inequality
 \begin{equation}\label{5}
f(y)\le (f(p^{-j_0}y))^{k_{j_0}}\le (f(p^{-j_0}y))^{p^{j_0}},\ \ y\in Y.
 \end{equation}
 Set $m=p^{j_0}$ and rewrite (\ref{5}) in the form
 \begin{equation}\label{22}
f(my)\le (f(y))^m,\ \ y\in Y.
 \end{equation}
Observe that we can assume that $m>2$ in  (\ref{22}). Indeed, if    $m=2$,
then $p=2$ and $j_0=1$. In this case   (\ref{21}) implies that $k_1>2$,
and in view of (\ref{4}), either $k_1=4$  or $k_1=3$. Let $k_1=4$.
Then (\ref{3}) takes the form $f(y)=(f(2^{-1}y))^4$. This implies that
$f(4y)=(f(y))^{16}$, and we may suppose that $m=4$ in (\ref{22}). Let
$k_1=3$. Then (\ref{3})   can be written as $f(y)=(f(2^{-1}y))^3\phi(y)$,
where $\phi(y)$ is  a characteristic function. It follows from this
that $f(4y)=(f(y))^9(\phi(2y))^3\phi(4y)$,  and we may also suppose
that $m=4$ in (\ref{22}).

 It may be noted that if $\mu=m_X$, then the theorem is proved.
Therefore, suppose that $\mu\ne m_X$. Then there is an element
$y_0\in Y$, $y_0\ne 0$,   such that $f(y_0)>0$.   Consider a group $F$ of the form (\ref{26_09_1}).
Since $Y$ is a group  with unique division by $p$, we can consider
on
the subgroup $F$ the function $g(r)=f(ry_0)$, $r\in F$.  Verify
that   (\ref{22}) implies the uniform continuity of the function
$g(r)$ in the topology induced on $F$ by the topology of  $\mathbb{R}$.
Taking into account that $g(r)$ is a positive definite function on
$F$,  for this it is enough to show that the function $g(r)$ is
continuous at zero.  We will follow the scheme of the proof of Lemma 10 in  \cite{TVP1991} and Lemma 8 in \cite{Fe11}.

Obviously,  it follows from (\ref{22})   that  $f(m^ky_0)\le (f(y_0))^{m^k}$
for all natural $k$. Hence,
$$
f(m^{-k}y_0)\ge (f(y_0))^{1/m^k}.
$$
It implies that
 \begin{equation}\label{21.1.1}
\textstyle g\left({1\over m^{k}}\right)\ge (g(1))^{1/m^k}.
 \end{equation}
Let $a>0$. Note that for all $x>0$ the inequality $1-e^{-ax}\le ax$ holds.  Taking this into account we get from   (\ref{21.1.1}) that
 \begin{equation}\label{7}
\textstyle 1-g\left({1\over m^{k}}\right)\le {C\over m^k},
 \end{equation}
 where $C=-\ln g(1)$.

 It is  obvious that for any arbitrary real-valued characteristic
function $f(y)$ the inequality
 \begin{equation}\label{31}
1-f(y_1+y_2)\le 2[(1-f(y_1))+(1-f(y_2))], \ \ y_1, y_2\in Y,
 \end{equation}
is fulfilled. By induction, this implies the inequality
 \begin{equation}\label{8}
1-f(y_1+\dots +y_q)\le \sum_{j=1}^q2^j(1-f(y_j)), \ \ y_j \in Y.
 \end{equation}
Let $i$ be an integer such that $0<i<m^k$. We have
\begin{equation}\label{9}
{i\over m^k}={i_1\over m}+\dots+{i_k\over m^k},
 \end{equation}
where $i_l$ are integers such that $0\le i_l< m$, $l=1, \dots, k$.  Substituting $q=i_l$,
$y_1=\dots=y_{i_l}={y_0\over m^{l}}$ in   (\ref{8}), we get
\begin{equation}\label{10}
\textstyle{1-g\left({i_l\over  m^{l}}\right)\le
\sum\limits_{j=1}^{i_l}2^j\left(1-g\left({1\over m^l}\right)\right)
\le 2^m\left(1-g\left({1\over m^l}\right)\right)}.
\end{equation}
Suppose that
\begin{equation}\label{101}
{i\over m^k}<{1\over m^s}.
\end{equation}
It follows from this that $i_l=0$ for all $l\le s$ in (\ref{9}). In this case, taking into
account (\ref{7}),  (\ref{9}) and  (\ref{10}),    we find from (\ref{8}) \begin{equation}\label{11}
\textstyle{1-g\left({i\over m^k}\right)=
1-g\left({i_{s+1}\over m^{s+1}}+\dots+{i_k\over m^k}\right)\le
\sum\limits_{j=1}^{k-s}2^j\left(1-g\left({{i_{s+j}\over m^{s+j}}}\right)
\right)}$$$$\textstyle{\le \sum\limits_{j=1}^{k-s}2^{j}2^m\left(1-g\left({{1\over m^{s+j}}}
\right)\right)
\le\sum\limits_{j=1}^{k-s}{C2^{m+j}\over m^{s+j}}=
{C2^{m}\over m^s}\sum\limits_{j=1}^{k-s}\left({2\over m}\right)^{j}}.
\end{equation}

In view of   $m>2$, we find from   (\ref{11})
that the inequality
\begin{equation}\label{12}
\textstyle 1-g\left({i\over m^k}\right)\le {C2^{m+1}\over (m-2)}{1\over m^s}
\end{equation}
is valid. Obviously, (\ref{101}) and (\ref{12}) imply the   continuity
at zero of the function  $g(r)$ in the topology induced on $F$ by
the topology of  $\mathbb{R}$, and hence,  the uniform continuity of   $g(r)$ in this topology.

It follows from (\ref{3}) that the function  $g(r)$ satisfies the equation
 \begin{equation}\label{23}
g(r)=(g(p^{-1}r))^{k_1}\cdots (g(p^{-l}r))^{k_l}, \ \ r\in F.
\end{equation}
We extend the function $g(r)$ by continuity from $F$ to a continuous
positive definite function on $\mathbb{R}$. We keep the notation $g$
for the extended function. It follows from (\ref{23}) that the function $g(s)$    satisfies
the equation
\begin{equation}\label{24}
g(s)=(g(p^{-1}s))^{k_1}\cdots (g(p^{-l}s))^{k_l}, \ \ s\in \mathbb{R}.
\end{equation}
By the generalized P\'olya  theorem, it follows from  (\ref{4})
and (\ref{24}) that $g(s)$  is the characteristic function of a symmetric
Gaussian distribution on the real line, and hence, we get the representation
\begin{equation}\label{32}
f(y)=\exp\{-\sigma(y_0) y^2\}, \ \ y=ry_0, \ \ r\in F,
\end{equation}
where $\sigma(y_0)\ge 0$.

When we obtained   (\ref{32}) we fixed an element  $y_0\in Y$, $y_0\ne 0$,
such that $f(y_0)>0$.
If we fix another element $\tilde y\in Y$, $\tilde y\ne 0$, such
that $f(\tilde y)>0$, we will obtain the
representation  $$f(y)=\exp\{-\sigma(\tilde y) y^2\}, \ \ y=r\tilde y,\ \
r\in F,$$
where $\sigma(\tilde y)\ge 0$. Since the subgroups
$\{y=ry_0:r\in F\}$ and  $\{y=r\tilde y:r\in F\}$ have a nonzero
intersection,  this implies that $\sigma(y_0)=\sigma(\tilde y)=\sigma$.
Put $E=\{y\in Y: f(y)\ne 0\}$. Thus, we have the representation
\begin{equation}\label{90}
f(y)=\exp\{-\sigma  y^2\}, \ \ y\in E.
\end{equation}

We will verify that $E$ is a subgroup of $Y$. First, observe that when we
got  (\ref{32}) we proved that if   $y\in E$, then $ry\in E$
for all $r\in F$.
Take $y_1, y_2\in E$. It follows from (\ref{90})
that   $f(p^{-k} y_j)\rightarrow 1$, $j=1, 2$,  as
$k\rightarrow\infty$, and then (\ref{31}) implies that
$p^{-k} (y_1+y_2)\in E$ for   sufficiently large natural
numbers $k$.  Hence  $y_1+y_2\in E$.

Thus, we obtain  the representation
\begin{equation}\label{33}
f(y)=
\begin{cases}
\exp\{-\sigma  y^2\}, & \text{\ if\ }\   y\in E,
\\  0, & \text{\ if\ }\ y\not\in
E.
\end{cases}
\end{equation}
Put $K=A(X, E)$. Let $\gamma$ be a Gaussian distribution on the group $X$
with the characteristic function  $\hat\gamma(y)=
\exp\{-\sigma  y^2\}$, $y\in Y$. Since $E=A(Y, K)$, it follows from (\ref{11a})
and (\ref{33}) that $f(y)=\hat\mu(y)=\hat\gamma(y)\hat m_K(y)$, and hence $\mu=\gamma*m_K$. Note that the subgroup $E$ has the property: if $py\in E$, then $y\in E$. This implies that $f_p(K)=K$.

 Get rid of the restriction that $\hat\mu(y)\ge 0$. Put $\nu=\mu*\bar\mu$.
Then $\hat\nu(y)=|\hat\mu(y)|^2\ge 0$ for all $y\in Y$, and the characteristic
function $\hat\nu(y)$ also satisfies equation  (\ref{2}). Then, as has
been proved above, the function   $\hat\nu(y)$ is represented in the
form (\ref{33}). On the one hand,  the subgroup $E$ possesses the property: if
$y\in E$, then  ${1\over p}y\in E$. Therefore   $E$ is not
isomorphic to the group $\mathbb{Z}$. Taking into account that if a subgroup of the group of rational numbers is not isomorphic to   $\mathbb{Z}$, then it is isomorphic to a group of the form (\ref{04.05.1}), this implies that   the group
$E$  is topologically isomorphic to a discrete additive
group $H_{\text{\boldmath $b$}}$ of the rational numbers of the form
$$
H_{\text{\boldmath $b$}}=
\left\{{m \over b_0b_1\cdots b_n} : \ n = 0, 1,\dots; \ m
\in {\mathbb{Z}} \right\}
$$
for some \text{\boldmath $b$}=$(b_0, b_1,b_2,\dots)$, where all  $b_j \in {\mathbb{Z}}$, $b_j > 1$. On the other hand, since $E=A(Y, K)$,
the group $E$ is topologically isomorphic to the character group of the
factor-group    $X/K$. By the Pontryagin duality theorem, the factor-group  $X/K$
is topologically isomorphic to the corresponding
{{\text{\boldmath $b$}}}-adic solenoid $\Sigma_{{\text{\boldmath $b$}}}$,
and hence the  group  $X/K$ contains no subgroup topologically isomorphic
to the circle group   $\mathbb{T}$.  Applying Lemma \ref{le1} to the
 group  $X/K$, and taking into account that each character of the subgroup
$E$ can be written as $(x, y)$, $x\in X$, $y\in E$,
we get the statement of the theorem in the general case from
representation (\ref{33}) for the function $\hat\nu(y)$. $\blacksquare$
 \begin{remark}\label{r1} Let $\Sigma_{{\text{\boldmath $a$}}}$ be
an   {{\text{\boldmath $a$}}}-adic solenoid. Let $\gamma$ be a
Gaussian distribution on   $\Sigma_{{\text{\boldmath $a$}}}$ with the characteristic function   of the form
\begin{equation}\label{12.02.1}
\hat\gamma(y)= \exp\{-\sigma y^2\}, \ \  y\in H_{\text{\boldmath $a$}},
\end{equation}
where   $\sigma\ge 0$.
Let $\alpha_j$, $j=1, 2, \dots, n$, $n\ge 2$, be topological automorphisms
of the group   $\Sigma_{{\text{\boldmath $a$}}}$ satisfying
the condition (\ref{1}). It is obvious that the characteristic function $\hat\gamma(y)$ satisfies equation (\ref{2}).

 Let $\Sigma_{{\text{\boldmath $a$}}}$ be an   {{\text{\boldmath $a$}}}-adic
solenoid with the property that there is a  prime number
$p$ such that    $\Sigma_{{\text{\boldmath $a$}}}$ contains no elements
of order  $p$. This implies  that
$\Sigma_{{\text{\boldmath $a$}}}$ is a group   with unique division by $p$.
Let $\alpha_j$, $j=1, 2, \dots, n$, $n\ge 2$,  be topological
automorphisms of the group $\Sigma_{{\text{\boldmath $a$}}}$ of the
form  $\alpha_j=\pm f^{-1}_{p^{k_j}}$, where $k_j$ are   natural numbers.
Let $K$ be a compact subgroup of  $\Sigma_{{\text{\boldmath $a$}}}$ such
that $f_p(K)=K$. Let $\xi_j$   be independent identically distributed
random variables with values in    $\Sigma_{{\text{\boldmath $a$}}}$
and distribution  $m_K$. Note that the condition $f_p(K)=K$ implies that
if $py\in A(H_{\text{\boldmath $a$}}, K)$, then $y\in A(H_{\text{\boldmath $a$}}, K)$. Taking this into account  and
using (\ref{11a}) it is not difficult to verify that the characteristic function $\hat m_K(y)$ also satisfies equation (\ref{2}).
 In view of Lemma \ref{le2}, it follows from has been said  the following statement.

{\it Let
$\Sigma_{{\text{\boldmath $a$}}}$ be an
{{\text{\boldmath $a$}}}-adic solenoid with the property that there is a prime number
$p$ such that    $\Sigma_{{\text{\boldmath $a$}}}$ contains no elements
of order  $p$. Let
$\alpha_j$, $j=1, 2, \dots, n$, $n\ge 2$, be topological
automorphisms of  $\Sigma_{{\text{\boldmath $a$}}}$ of the form $\alpha_j=\pm f^{-1}_{p^{k_j}}$, where $k_j$ are   natural numbers. Assume that $\alpha_j$ satisfy  the condition
$(\ref{1})$. Let $\gamma$ be a    Gaussian distribution on
$\Sigma_{{\text{\boldmath $a$}}}$ with the characteristic function of the form $(\ref{12.02.1})$, and $K$ be a compact subgroup of
$\Sigma_{{\text{\boldmath $a$}}}$ such that $f_p(K)=K$. Put
$\mu=\gamma*m_K$. Let $\xi_j$ be independent identically distributed
random variables with values in    $\Sigma_{{\text{\boldmath $a$}}}$
and distribution        $\mu$.  Then   $\xi_j$  and
$ \alpha_1\xi_1+\dots+ \alpha_n\xi_n$ are identically distributed.}

Hence,
we cannot narrow down  the class of distributions in Theorem \ref{th1}
which is characterized by  the property of equidistribution
of    $\xi_j$ and     $\alpha_1\xi_1+\dots+\alpha_n\xi_n$.
\end{remark}
It turns out that Theorem \ref{th1} is false if the condition $(i)$
is not satisfied. Namely, the following statement holds.
\begin{proposition}\label{pr1} Let $X=\Sigma_{{\text{\boldmath $a$}}}$
be an {{\text{\boldmath $a$}}}-adic solenoid satisfying the condition:

$(i)$ There are   two   prime numbers $p$ and $q$
such that the group $X$ contains no elements of order   $p$ and $q$.

Then there are topological automorphisms $\alpha_j$,
$j=1, 2, \dots, n$, $n\ge 2$, of the group   $X$ satisfying the condition
$(\ref{1})$, and independent identically distributed random
variables $\xi_j$ with values in    $X$ and distribution  $\mu$
such that  $\xi_j$ and    $\alpha_1\xi_1+\dots+\alpha_n\xi_n$
are identically distributed, whereas  $\mu$ can not be represented
as a convolution  $\mu=\gamma*m_K$, where $\gamma\in \Gamma(X)$,
and $K$ is a compact subgroup of the group   $X$.
\end{proposition}
Theorem  \ref{th1} and Proposition \ref{pr1} imply the following statement.
\begin{corollary}\label{co2} Let $X=\Sigma_{{\text{\boldmath $a$}}}$
be an {{\text{\boldmath $a$}}}-adic solenoid. Let $\alpha_j$, $j=1, 2, \dots, n$, $n\ge 2$, be   topological automorphisms of the group   $X$  satisfying
the condition $(\ref{1})$.
Let $\xi_j$ be independent identically distributed random variables
with values in    $X$ and distribution  $\mu$. The equidistribution of  $\xi_j$ and     $\alpha_1\xi_1+\dots+\alpha_n\xi_n$ implies  that
$\mu=\gamma*m_K$, where $\gamma\in \Gamma(X)$
and $K$ is a compact subgroup of the group   $X$, if and only if the  condition
$(i)$ of Theorem $\ref{th1}$ is satisfied.
\end{corollary}
To prove Proposition \ref{pr1} we need the following lemma
(\!\!\cite[(32.43)]{HeRo2}).
\begin{lemma}\label{le3} Let $Y$
be an   Abelian group, let $L$
be a subgroup of $Y$, and let $g(y)$ be a positive definite function on
$L$. If a function
$f(y)$ on the group $Y$ is represented in the form
 $$
 f(y)=
 \begin{cases}
 g(y), & \text{\ if\ }\ \ y\in L,
 \\ 0, & \text{\ if\ }\ \ y\not\in L,
 \end{cases}
 $$
then $f(y)$ is a   positive definite function.
\end{lemma}
{\bf Proof of Proposition \ref{pr1}.}  Suppose for definiteness that $p<q$.
In order not to complicate
the notation, we will identify $Y$ with $H_{\text{\boldmath $a$}}$.
The condition $(i)$ in Proposition \ref{pr1} implies
that $f_p, f_q\in {\rm Aut}(X)$, and hence,  $f_p, f_q\in {\rm Aut}(Y)$.
Therefore both $X$ and $Y$ are groups with unique division by  $p$ and by $q$.
Denote by $H$ a subgroup of $Y$  of the form
$$H= \left\{{m \over n}\in Y: n  \  \text{is not divided by} \  p;   \ m
\in {\mathbb{Z}} \right\}.
$$
Let $k$ be a natural number. Set $$H_k=\left\{{m\over p^kn}\in Y: m, n  \
\text{are not divided by} \ p\right\}$$
and
$$   L=H\cup H_1.$$
It is obvious that $Y=H\cup\bigcup\limits_{k=1}^\infty H_k$ and $L$ is a
subgroup of   $Y$.  Denote by  $G$ the character group of
the group   $L$, and  put $F=A(G, H)$. It is easy to see
that  $L/H\cong \mathbb{Z}(p)$, where $\mathbb{Z}(p)$ is the group of
residue  classes modulo  $p$. Since the character group of the
factor-group    $L/H$ is topologically isomorphic to the annihilator
$A(G, H)$, we have $F\cong \mathbb{Z}(p)$. Take $0<c<1$. Let $\omega$ be a distribution
on the group   $F$ of the form $\omega=cE_0+(1-c)m_F$, where $E_0$ is the
distribution concentrated at zero.  Consider  $\omega$ as a distribution
on the group     $G$. Then     (\ref{11a}) implies that the characteristic
function  $\hat\omega(l)$, $l\in L$,  is of the form
$$\hat\omega(l)=
 \begin{cases}
 1, & \text{\ if\ }\ \ l\in H,
 \\ c, & \text{\ if\ }\ \ l\in H_1.
   \end{cases}
 $$
 Obviously, $\hat\omega(l)$  is a positive definite function. Consider
on the group
$Y$ the function
$$f(y)=
 \begin{cases}
 \hat\omega(y), & \text{\ if\ }\ \ y\in L,
   \\ 0, & \text{\ if\ }\ \ y\not\in L.
 \end{cases}
 $$
  By Lemma \ref{le3}, $f(y)$ is a positive definite function. By the
Bochner theorem, there is a probability distribution   $\mu$ on the group $X$
such that $f(y)=\hat\mu(y)$. It is obvious that the distribution  $\mu$ is
not  represented as a convolution  $\mu=\gamma*m_K$,
where $\gamma\in \Gamma(X)$, and $K$ is a compact subgroup of
  $X$.

Take   natural numbers   $l$ and $m$, where $l< m$, such that the remainders
of the division of the numbers $q^{2l}$ and $q^{2m}$ by $p^2$ are equal.
Then $q^{2m}-q^{2l}$ is divided by   $p^2$. Hence, $q^{2m-2l}-1$ is also
divided by $p^2$, and the remainder of the division  $q^{2m-2l}$ by $p^2$ is
equal to 1. Put $a=m-l$. We have
\begin{equation}\label{71}
q^{2a}=p^2b+1.
\end{equation}
Put $\alpha_1=\dots=\alpha_b=f_pf^{-1}_{q^a}$,
$\alpha_{b+1}=f^{-1}_{q^a}$, $n=b+1$. In view of  (\ref{71}),  the
topological automorphisms   $\alpha_j$,  $j=1, 2, \dots, n$,    of the
group $X$  satisfy the condition $(\ref{1})$.

Let $\xi_j$, $j=1, \dots n$, be independent identically distributed
random variables with values in    $X$ and distribution  $\mu$.
Verify that  the function $f(y)$ satisfies the equation
  \begin{equation}\label{51}
\textstyle f(y)=\left(f\left({p\over q^a}y\right)\right)^bf\left({1\over q^a}y\right),  \ \ y\in Y.
\end{equation}
Indeed, if $y\in H$, then ${p\over q^a}y, {1\over q^a}y\in H,$ and both
parts of equation  (\ref{51}) are equal to 1. If $y\in H_1$,
then ${p\over q^a}y\in H,$ ${1\over q^a}y\in H_1,$ and both parts of
equation  (\ref{51}) are equal to $c$. If $y\in H_k$,  $k\ge 2$,
then ${1\over q^a}y\in H_k,$ and both parts of equation  (\ref{51}) are
equal to zero. Thus, the function $f(y)$ satisfies   equation (\ref{51}),
and by Lemma \ref{le2},   $\xi_j$ and
$\alpha_1\xi_1+\dots+\alpha_n\xi_n$ are identically distributed. $\blacksquare$
\begin{remark}\label{r3} Let $X$ be a locally compact Abelian group and
let $\nu\in{\rm  M}^1(X)$. Put
$B=\{y\in Y:\hat\nu(y)=1\}$. It is well-known that  $B$
is a closed subgroup of $Y$ and  the distribution $\nu$ is supported
in the subgroup $A(X, B)$.  Therefore, in the notation of Proposition
\ref{pr1}   the distribution    $\mu$ constructed in the proof of
Proposition \ref{pr1} is supported in the annihilator $A(X, H)$, i.e in   a   proper closed subgroup of the group $X$.

We can strengthen Proposition \ref{pr1} as follows. Let $X=\Sigma_{{\text{\boldmath $a$}}}$
be an {{\text{\boldmath $a$}}}-adic solenoid satisfying the condition
$(i)$ of Proposition \ref{pr1}.  Let $\gamma$ be
a non-degenerate   Gaussian distribution on  $X$ with the
characteristic function of the form (\ref{12.02.1}).
Set $\lambda=\gamma*\mu$, where the distribution    $\mu$ constructed in the proof of
Proposition \ref{pr1}. Then the characteristic function of the
distribution $\lambda$ is of the form
$$\hat\lambda(y)=
 \begin{cases}
 \exp\{-\sigma y^2\}, & \text{\ if\ }\ \ y\in H,
   \\ c\exp\{-\sigma y^2\}, & \text{\ if\ }\ \ y\in H_1,
   \\ 0, & \text{\ if\ }\ \ y\not\in L.
 \end{cases}
 $$
It is obvious that the distribution $\lambda$  can not be represented
as a convolution of a Gaussian distribution on $X$ and the Haar
distribution of a compact subgroup of the group   $X$. Let $\xi_j$
be independent identically distributed random variables with values
in    $X$ and distribution   $\lambda$. Arguing   as in the proof
of Proposition \ref{pr1} and using Lemma \ref{le2}, we are convinced
that   $\xi_j$ and   $\alpha_1\xi_1+\dots+ \alpha_n\xi_n$
are identically distributed. As is easily seen, the distribution $\lambda$ is not supported
in any  proper close subgroup of the group $X$.
\end{remark}
\begin{remark}\label{r4}
On the real line, both the Skitovich--Darmois theorem and the generalized P\'olya theorem characterize the Gaussian distribution. Corollary \ref{co2} gives a complete description of {{\text{\boldmath $a$}}}-adic solenoids on which the generalized P\'olya theorem is valid. It is interesting to note that, as follows from   \cite {FG1}, there are no
{{\text{\boldmath $a$}}}-adic solenoids on which the Skitovich--Darmois theorem holds. Namely, the following statement takes place.

{\it Let $X=\Sigma_{{\text{\boldmath $a$}}}$ be an   {{\text{\boldmath $a$}}}-adic solenoid. Then there are topological automorphisms   $\alpha_j, \beta_j$
of the group $X$ and
independent random variables $\xi_1$ and $\xi_2$ with values in
$X$ and distributions $\mu_1$ and   $\mu_2$ such that the linear forms $L_1 = \alpha_1\xi_1 +
\alpha_2\xi_2$ and $L_2 = \beta_1\xi_1 + \beta_2\xi_2$ are independent, whereas   $\mu_j$ can not be represented
as   convolutions  $\mu_j=\gamma_j*m_{K_j}$, where  $\gamma_j\in \Gamma(X)$, and $K_j$ are compact subgroups of the group  $X$, $j=1, 2$.}
\end{remark}
\begin{remark}\label{r5}
As noted earlier each connected locally compact Abelian group of dimension 1 is topologically isomorphic to either the additive group of real numbers, or the circle group, or an \text{\boldmath $a$}-adic
solenoid.   The generalized P\'olya theorem  on the additive group of real numbers  characterizes   Gaussian distributions (Theorem A). The generalized P\'olya theorem on {{\text{\boldmath $a$}}}-adic solenoids $\Sigma_{{\text{\boldmath $a$}}}$ satisfying the condition $(i)$ of Theorem \ref{th1} characterizes convolutions  of Gaussian distributions on $\Sigma_{{\text{\boldmath $a$}}}$ and Haar distributions on compact subgroups of $\Sigma_{{\text{\boldmath $a$}}}$  (Theorem \ref{th1}). Discuss the generalized P\'olya theorem on the circle group.

Let $\alpha_j$, $j=1, 2, \dots, n$, $n\ge 2$, be topological automorphisms
of the circle group $\mathbb{T}$, and let $\xi_j$ be independent identically distributed random variables with
values in    $\mathbb{T}$ and distribution  $\mu$. Consider the linear form $L=\alpha_1\xi_1+\dots+\alpha_n\xi_n$.
Taking into account that $I$ and $-I$ are the only   topological automorphisms of $\mathbb{T}$,    $L$ is of the form  $L=\pm\xi_1+\dots+\pm\xi_n$. Assume for definiteness that $L=\xi_1+\dots+\xi_m-\xi_{m+1}-\dots-\xi_n$. Set $f(y)=\hat\mu(y)$. If $\xi_j$ and $L$ are identically distributed, then  by Lemma \ref{le2}, the function $f(y)$ satisfies equation  (\ref{2}) which takes the form
\begin{equation}\label{06.05.1}
f(y)=(f(y))^m(f(-y))^{n-m},
\ \ y\in \mathbb{Z}.
 \end{equation}
 Put $E=\{y\in \mathbb{Z}: f(y)\ne 0\}$. It follows from (\ref{06.05.1}) that $|f(y)|=1$ for all $y\in E$. This implies that $E$ is a subgroup of $\mathbb{Z}$ and there is an elements $x\in \mathbb{T}$ such that $f(y)=(x, y)$  for all $y\in E$.   Thus, we have
 \begin{equation}\label{06.05.2}
 f(y)=
 \begin{cases}
 (x, y), & \text{\ if\ }\ \ y\in E,
   \\ 0, & \text{\ if\ }\ \ y\not\in E.
 \end{cases}
  \end{equation}
Put $K=A(\mathbb{T}, E)$. Then $E=A(\mathbb{Z}, K)$. Taking into account (\ref{11a}), it follows from  (\ref{06.05.2}) that $\mu$ is a shift by the element $x$ of the Haar distribution $m_K$. We see that the generalized P\'olya theorem on the circle group $\mathbb{T}$ characterizes shifts of Haar distributions on compact subgroups of $\mathbb{T}$.

Thus, we have completely studied the generalized P\'olya theorem on connected locally compact Abelian groups of dimension 1.

\end{remark}

\vskip 1 cm

\noindent B. Verkin Institute for Low Temperature Physics and Engineering \\ of the National Academy of Sciences of Ukraine

\bigskip

47, Nauky ave, Kharkiv, 61103, Ukraine \\

\bigskip

E-mail: feldman@ilt.kharkov.ua

\bigskip


\newpage

\begin{thebibliography}{99}

\bibitem{Fe1}    Feldman, G. M.:{  Gaussian distributions on locally compact Abelian groups}.   Theory
Probab. Appl.  {\bf 23},    529--542 (1979)

\bibitem{TVP1991}  Feldman, G. M.:   Cauchy distribution on Abelian
groups and its characterization.   Theory
Probab. Appl.  {\bf 36},    670--681  (1991)

\bibitem{Fe9} Feldman, G. M.: Arithmetic of probability
    distributions and characterization  problems on Abelian
    groups. Transl. Math. Monographs
{\bf 116},     Amer. Math. Soc., Providence, RI  (1993)


\bibitem{Fe11} Feldman, G. M.:  On a characterization theorem for
the Cauchy distribution on abelian groups.    Siberian Math. J.  {\bf 36},  196--201    (1995)

\bibitem{Fe-SD-2003}  Feldman, G. M.: {A characterization of the Gaussian
distribution on Abelian groups}. Probab. Theory Relat. Fields
 {\bf 126},     91--102 (2003)

\bibitem{Fe4} Feldman, G. M.: { On a characterization theorem for locally
compact abelian groups}.   Probab. Theory Relat. Fields
 {\bf 133},    345--357 (2005)

 \bibitem{Fe5} Feldman, G. M.: Functional equations and characterization problems
on locally compact Abelian groups. EMS Tracts in Mathematics,  {\bf
5}  European Mathematical Society, Zurich  (2008)

\bibitem{Fe20bb} Feldman, G. M.: The Heyde theorem for  locally compact Abelian
groups.  J.   Funct. Anal.   {\bf 258},  3977--3987 (2010)

\bibitem{Fe6} Feldman, G. M.: { On a characterization of convolutions of  Gaussian and Haar distributions}.
Math. Nachr.  {\bf 286},  340--348  (2013)

\bibitem{JOTP2015} Feldman, G. M.: {On the Skitovich--Darmois theorem for the group of p-adic numbers}. J. of Theor. Probab.   {\bf 28},  539--549  (2015)

\bibitem{Solenoid} Feldman, G. M.: On a characterization theorem for
connected locally compact Abelian groups. J. Fourier Anal. Appl.   {\bf 26}  (14),  1--22  (2020)

\bibitem{FG1} Feldman, G. M., Graczyk, P.:      On the Skitovich-Darmois theorem
on compact Abelian groups.  J. of Theor. Probab.   {\bf 13}, 859--869  (2000)

 \bibitem{FG} Feldman, G. M., Graczyk, P.: {The Skitovich-Darmois theorem for locally compact Abelian groups}. J. of the Australian Mathematical Society   {\bf 88},     339--352  (2010)

\bibitem{MiFe5} Feldman, G. M.,  Myronyuk, M.V.: {Independent linear forms on
connected Abelian groups}. Math. Nach.   {\bf 284},      255--265 (2011)

\bibitem{HeRo1} {Hewitt, E.,  Ross, K.A.:} Abstract
   Harmonic Analysis.   {\bf 1}.   Grundlehren Math. Wiss.   {\bf 115},
Berlin --
    Gottingen -- Heildelberg (1963)

\bibitem{HeRo2}  {Hewitt, E.,  Ross, K.A.:} Abstract
   Harmonic Analysis.   {\bf   2}.   Grundlehren Math. Wiss, {\bf  152}.
Springer--Verlag,   New York -- Berlin  (1970)

\bibitem{KaLiRa}    Kagan, A. M.,  Linnik, Yu. V.,  Rao, C.R.:
{  Characterization problems in mathematical statistics}. Wiley
Series in Probability and Mathematical
Statistics, John Wiley $\&$ Sons, New York-London-Sydney (1973)

\bibitem{M} Myronyuk, M.V.:  Independent linear forms on the group $\Omega_p$. J. of Theor. Probab.   {\bf 33},     1--21    (2020)

\bibitem{Pa} Parthasarathy, K.R.:  {Probability measures
on metric spaces}. Academic Press,  New York and London  (1967)

\end{thebibliography}
\end{document}